\documentclass[12pt]{article}
\usepackage{amsmath, amssymb,paralist,amsthm}

\def\E{\mathbb{E}}
\def\R{\mathbb{R}}
\def\N{\mathbb{N}}

\newtheorem{Theorem}{Theorem}[section]

\newtheorem{Remark}[Theorem]{Remark}

\newtheorem{Lemma}[Theorem]{Lemma}

\numberwithin{equation}{section}

\begin{document}

\centerline{\Large \bf The turnpike theorems for Markov games \footnote{Supported by the AFOSR grant FA9550-09-1-0664 'Nonlinear Markov control processes and games'} \footnote{Submitted to { \it Dynamic Games and Applications} for publication}}

\bigskip
\bigskip
\centerline{\bf Vassili Kolokoltsov and Wei Yang}
\smallskip
\centerline{\small Department of Statistics, University of Warwick}
\centerline{\small  Coventry,  CV4 7AL,   UK}
\centerline{\small  \it v.kolokoltsov@warwick.ac.uk \quad wei.yang@warwick.ac.uk}

\bigskip

\begin{abstract}
This paper has a two-folded purpose. First, we attempt to outline the development of the turnpike theorems in the the last several decades. Second, we study turnpike theorems in  finite-horizon two-person zero-sum Markov games on a general Borel state space. Utilising the Bellman (or Shapley) operator defined for this game, we prove the stochastic versions of the early turnpike theorem on the set of optimal strategies and the middle turnpike theorem on the distribution of the state space.
\end{abstract}

{\medskip\par\noindent

}

\section{Introduction}

Turnpike theorems are a set of economic theorems on optimal control problems, which characterise solutions (optimal paths) of optimal control problems and focus on the stationary optimal one. The {\it turnpike theorem} was  first proposed  by Dorfman, Samuelson and Solow in \cite{DSS1958} within the context of an optimal growth model. In the literature, there are extensive work on the turnpike theorems in many diversified settings of optimal control problems and games, ranging from discrete time to continuous time, from the deterministic case to the stochastic case, with or without discounting. There is no clear-cut and unified classification for turnpike theorems. For example, McKenzie \cite{Mc1976} distinguished three kinds of turnpikes: the early, late and middle turnpikes, but in \cite{M1998}, he used a different classification scheme with the {\it Samuelson} and {\it Ramsey turnpike}. We think that this is helpful to clarify different types of turnpike theorems and to give a clear picture on the state of this subject. Hence we start this paper with a survey on the main developments of the turnpike theorems in the last several decades.

Further on, in this paper, we study turnpike theorems in finite-horizon  zero-sum Markov games on a general Borel state space following the approach developed in \cite{K1992} in which a finite state space was considered. A general approach to Borel state space stochastic games was developed by Maitra and Sudderth \cite{MS1993}.  They proved the existence of the game value for general universally measurable strategies. However, the optimal strategies need not exist in their games and it is assumed the transition probability function is strongly continuous in the actions of one player.

In our setting, under the weaker assumption \eqref{prop2} in Section \ref{S3},  we show the convergence of the long time average cost (reward) of the game and prove the turnpike theorems on the set of strategies and on the distribution of the states.

The organisation of this paper is as follows. In Section \ref{survey}, we sketch a picture of turnpike theory showing how the theory was developed. In Section \ref{modelling}, we present the basic setting for a two-person zero-sum Markov game and define the Bellman (or Shapley) operator. The main results will be presented and proved in Section \ref{S3}.  Section \ref{conclusion} concludes this paper and proposes further research.

\section{Recent development in the turnpike theory}\label{survey}

In this section, we outline the development of the turnpike theory in the last several decades.  We shall start with turnpike theorems in the deterministic setting, followed by the stochastic analog of the turnpike theorems.

\subsection{The turnpike theorems in the deterministic setting}

For the turnpike theorems in the deterministic setting, it is important to distinguish the {\it asymptotic} and {\it neighborhood turnpike theorems} for solutions to {\it discounted infinite-horizon} optimal control problems from the {\it classical turnpike theorems} for solutions to {\it undiscounted finite-horizon} problems. Some surveys in the deterministic setting can be also found in \cite{M1998, Z2006, B2007, KP2011}.

The {\it classical turnpike theorem} considers a deterministic finite-horizon ($T>0$) optimal control problem on a given state space $X$. In a discrete time case, the problem is
\begin{gather}\label{DWD}
\begin{cases}
\text{maximise}\quad&\sum_{t=0}^{T-1}  L(x_t, x_{t+1}) \\[.4em]
\text{subject to} \quad &x_0, x_T \in X\, \text{given},\quad \{ x_t\}_{t=0}^T\subset X
\end{cases}
\end{gather}
with a cost (utility) function $L:X^2\to R$, $T\in \N$.  In a continuous time case, the problem is
\begin{gather}\label{CWD}
\begin{cases}
\text{maximise}\quad&\int_0^T L(x_t, u_t)dt \\[.4em]
\text{subject to} \quad &\dot x_t=h(x_t, u_t),\quad u_t\in U\\
    &x_0, x_T \in X\, \text{given},\quad \{ x_t\}_{t\in [0,T]}\subset X
\end{cases}
\end{gather}
with a cost (utility) function $L:X\times U\to R$, where $h:X\times U\to R$ is a continuous mapping with respect to both arguments, $U$ is the set of  admissible controls and $T\in\R^+$. The initial value $x_0$ and the terminal value $x_T$ are held fixed. Let $x^*$ denote a solution of the maximization problem
$$\text{maximise}\quad L(x,x),\quad x\in X.$$
The path $x^*:=\{x_t: x_t=x^*, \forall t\in [0,T]\}$ is called a stationary optimal path.
Often, the uniqueness of the stationary optimal path is assumed.

In the case where the terminal value $x_T$  is not given, an additional term, interpreted as the terminal payoff function, should be added to the above problems. That is, the optimal problem is, in a discrete time setting
\begin{gather}\label{DWD2}
\begin{cases}
\text{maximise}\quad&\sum_{t=0}^{T-1}  L(x_t, x_{t+1}) +S(x_T)\\[.4em]
\text{subject to} \quad &x_0\, \text{given},\quad \{ x_t\}_{t=0}^T\subset X
\end{cases}
\end{gather}
and, in a continuous time case
\begin{gather}\label{CWD2}
\begin{cases}
\text{maximise}\quad&\int_0^T L(x_t, u_t)dt +S(x_T)\\[.4em]
\text{subject to} \quad &\dot x_t=h(x_t, u_t),\quad u_t\in U\\
    &x_0\, \text{given},\quad \{ x_t\}_{t\in [0,T]}\subset X
\end{cases}
\end{gather}
where the terminal payoff function $S:X\to \R$.

The original formulation of the continuous time problem \eqref{CWD} is ascribed to Samuelson-Solow \cite{SS1956} and the discrete time problem \eqref{DWD} to Gale \cite{G1967}  and Mckenzie \cite{M1968, Mc1976}.  It was further investigated by e.g.  Radner \cite{R1961}, Morishima \cite{M1961}, McKenzie  \cite{Mc1963} and Makarov and Rubinov \cite{MR1973},  Atsumi \cite{A1965} and Haurie \cite{H1976}. Mckenzie \cite{M1998} discussed in detail the development of the turnpike theory in both the von Neumann growth model and the Ramsey's growth model and specified a classification of Samuelson turnpike and Ramsey turnpike. The classical techniques which were developed to prove the classical turnpike theorems are critically dependent on the convexity of the control set $U$ and technology function $h$ and the preference function $L$. These classical techniques fall into two classes: Value-Loss methods and Monotonicity method, see more discussion in \cite{J1997}.

The classical turnpike theorem states that  if $T$ is large enough, the optimal path $\{x^*_t\}$ that transfers the system from $x_0$ to $x_T$ approaches the unique optimal stationary level $x^*$, and stays close to it for a large fraction of $T$, and moves away toward the terminal state only in the final periods; that is, for each $\epsilon>0$, there exist an integer $T_0$ such that for each $T\geq 2T_0$,
\begin{equation}\label{T1}
||x_t^*-x^*||\leq \epsilon, \quad \forall t\in [T_0, T-T_0].
\end{equation}
This is called the {\it middle turnpike theorem} in \cite{Mc1976}.  Moreover, if this result holds for any initial value $x_0\in X$, it is the {\it global} version; if the result holds for some initial value $x_0$ which is near enough to $x^*$, it is  the {\it local} version, e.g. \cite{BC2007}.  Further, if for each $\epsilon>0$, there exists an integer $T_0$ and a number $\delta>0$ such that $||x_0-x^*||\leq \delta$ implies \eqref{T1} for all $t\in [0, T-T_0]$ with $T\geq 2T_0$,  this was called the {\it early turnpike theorem} in \cite{Mc1976}. Since this image resembles a map of an intercity highway or turnpike with entrance and exit ramps, a stationary optimal path $x^*$ is often referred to as a {\it turnpike}.

A version of the classical turnpike theorem for games was also presented. For example, two-player zero-sum differential games were studied by Zaslavski \cite{Z1999}, showing that there exists a pair of optimal stationary paths $x^*$ and $y^*$ for the two players and presenting a version of classical middle turnpike theorem, and by Alvarez and Bardi \cite{AB2010,AB2007,B2009}, proving a turnpike theorem for the value of the game.

In the classical setting, the state space $X$ was often a compact, convex, bounded and closed subset of a finite-dimensional space; the cost (utility) function  $L$ is often  continuous, smooth and strictly concave (or convex); hence the turnpike theorem was proved by the concavity analysis. Later, under upper (lower) semicontinuity, some extensions to general Banach state spaces were obtained by e.g. Kolokoltsov \cite{K1992},  Yakovenko and L.A. Kontorer \cite{YK1992}, Zaslavski \cite{Z1995-2, Z1998_2, Z2004}, Mamedov \cite{M2002}.


When the future is discounted, the deterministic optimal control problem with infinite-horizon is described, in a discrete-time setting as
\begin{gather}\label{dis-disc}
\begin{cases}
 \text{maximise}\quad&\sum_{t=0}^\infty\delta^t L(x_t,x_{t+1})\\[.4em]
 \text{subject to}\quad &x_0\in X, \quad \{x_t\}_{t=0}^\infty\subset X
 \end{cases}
 \end{gather}
and in a continuous-time setting as
\begin{gather}\label{cont-disc}
\begin{cases}
 \text{maximise}\quad &\int_0^\infty e^{-\rho t}L(x_t, u_t)dt\\[.4em]
 \text{subject to}\quad &\dot x_t=h(x_t, u_t),\quad u_t\in U \\
            &x_0\in X, \quad \{x_t\}_{t\geq 0}\subset X
 \end{cases}
 \end{gather}
with a discount factor $0<\delta\leq 1$ and a discount rate $\rho\geq 0$. The case where $\delta=1$ and $\rho=0$ is referred to as the undiscounted version of the problem \eqref{dis-disc} and \eqref{cont-disc}. The asymptotic behaviour of the optimal path in the discounted framework depends critically on the magnitude of the initial stock and the discount factor.

Two classes of turnpike theorems should be distinguished when the future is discounted. The first class is the {\it asymptotic turnpike theorem}, see e.g. in Scheinkman \cite{S1976}, Bewley \cite{B1982} , Yano \cite{Y1984}, Carlson\cite{CJH1987}, Carlson, Haurie and Jabrane \cite{C1990}, Zaslavski \cite{Z1995-1}, Park  \cite{P2000}. It states that there exists a discount factor $\delta'\in(0,1)$ such that for $x_0\in X$ and any $\delta\in[\delta', 1)$, the optimal path $\{x^*_t\}$ starting at $x_0$ converges to the stationary path $x^*$.

The asymptotic turnpike theorem was proved under strictly concavity assumptions on $L$ by e.g. Scheinkman \cite{S1976} and Bewley \cite{B1982},  Majumdar, Nermuth \cite{MN1982}, Kamihigashi and Roy \cite{KR2006}. Majumdar and Mitra \cite{MM1982} interpreted the function $L$ as the consumption and set $L(x_t,x_{t+1}):=f(x_t)-x_{t+1}>0$ for all $t$ and considered a nonconvex (convex-concave) technology function $f$, i.e. $f$ is assumed to be strictly increasing and twice continuously differentiable such that $f$ satisfies  $f'(\infty)<1<f'(0)<\infty$ and there is a real number $k_1>0$ such that
$f''(x)=0$ for $x=k_1$, $f''(x)>0$ for $0\leq x<k_1$ and $f''(x)<0$ for $x>k_1$.

Montrucchio \cite{M1995, M1995-2} assumed certain curvature restrictions, i.e. the function $L$ is continuous and strongly $(\alpha, \beta)$-concave, i.e. $L(x,y)+\frac{1}{2}\alpha|x|^2+\frac{1}{2}\beta|y|^2$ is concave with $\alpha+\beta>0$  and that the value function $W(x_0):=\max\sum_{t=0}^\infty\delta^t L(x_t,x_{t+1})$ is concave-$\gamma$ for all $\delta \in (\delta_0,1)$, i.e. $W(x)+\frac{1}{2}\gamma|x|^2$ is convex for all $\delta \in (\delta_0,1)$. He proved that there exists $\delta'=1-\frac{\alpha+\beta}{\gamma}$ such that for a discount factor $\delta \in (\max [\delta_0,1-\frac{\alpha+\beta}{\gamma}],1) $ the local asymptotic turnpike theorem holds, where $\alpha+\beta$ is a measure of the lower curvature of $L$ and $\gamma$ is a measure of its upper curvature. Carlso and Haurie \cite{CH1995} proved a version of the asymptotic turnpike theorem for infinite horizon open-loop differential games with  $\delta=1$ under a ``strict diagonal concavity condition" on $L$. An {\it uniform asymptotic turnpike theorem} was proved by Cuong LeVan and Lisa Morhaim \cite{LVM2006} under similar assumptions, i.e. if the function $L(x,y)$ is a strictly concave function, increasing  in $x$ and decreasing in $y$, then for any initial point $x_0\in X$,  there exists a $\beta'\in(0,1)$ such that for any $\beta\in [\beta', 1)$, the asymptotic turnpike property holds.

Araujo and Scheinkman \cite{AS1977}, McKenzie \cite{Mc1977}, Dasgupta and Mckenzie \cite{DMc1985} and Hyun Park \cite{P2000} applied the {\it implicit function theorem} to prove the asymptotic turnpike theorem by assuming a dominant diagonal blocks condition:    for an infinite matrix $D$ formed by $n\times n$ blocks $D_{ij}$, $i,j\in \{1,2,\dots\}$, $D_{ii}$ is invertible, $\sup_i|D_{ii}^{-1}|<\infty$ and $\sup_i\sum_{j=1, i\neq j}^\infty |D_{ii}^{-1}D_{ij}|<1$, where the norm $|D_{ij}|=\sup_{|z|=1}|D_{ij}z|$.


The second class is the {\it neighborhood turnpike theorem}, see e.g. \cite{c1966, M1979, Y1984, C1985, Mc1986, FK1990, MM1999, KR2006}. It states that for any $\epsilon >0$, there exists $T'>0$ and $\delta'\in(0,1)$ such that, an optimal path stating from $x_0$ at the discount factor $\delta\in (\delta',1)$ eventually stays within the $\epsilon$-neighborhood of a stationary path $x^*$,  i.e.
$$||x^*_t-x^*||\leq \epsilon,\quad \forall t\geq T'.$$
If a smaller neighborhood $\epsilon$ is chosen, the closer the discount factor $\delta$ is to $1$.

To prove the neighborhood turnpike theorem, e.g. in Mckenzie \cite{Mc1986}, Yano \cite{Y1984}, instead of strictly concavity assumption, the differentiability of the function $L$ is required. A version of {\it uniform neighborhood turnpike theorem} was proved in Guerrero-Luchtenberg and Leonardo \cite{LL2000}. It states that, under certain differentiability and concavity conditions on $L$, the neighborhood turnpike theorem  holds for almost all initial states.

Finite-horizon problem with discounting was examined in e.g. \cite{FO1985, FK1990}. Based on convexity condition and using Implicit Programming formulation, Feinstein and Oren \cite{FO1985} proved a {\it funnel asymptotic turnpike theorem}, namely, an optimal path lies in an exponentially bounded region (determined by the discount factor) of the optimal stationary path. Fershtman and Kamien \cite{FK1990} considered a finite-horizon differential games with quadratic cost function and proved a neighborhood turnpike theorem for the feedback equilibrium strategies.

\subsection{The stochastic version of turnpike theorems}

The turnpike theorems under uncertainty was first developed by William and Mirman \cite{BM1972}. The uncertainty or the stochastic environment is represented by a filtered probability space $(\Omega, \mathcal{F},\nu; \{\mathcal{F}_t\}_{t\geq 0})$, $\omega\in \Omega$ is a possible state of the environment, $\mathcal{F}$ is the Borel $\sigma$-field, $\nu$ is a probability measure on $\Omega$, $\{\mathcal{F}_t\}_{t\geq 0}$ is the filtration satisfying the ``usual" conditions, i.e. it is right-continuous and complete. The underlying evolution of the states $\{x_t\}$ is a stochastic process on $(\Omega, \mathcal{F},\nu)$. Then the stochastic optimal control problem, in the undiscounted case, can be described, in the discrete-time setting as
\begin{equation}\label{DS_1}
\text{maximise}\quad \E\left[ \sum_{t=0}^T  L(x_t,u_t)\right]
\end{equation}
and in the continuous-time setting as
\begin{equation}\label{Ctn_1}
\text{maximise}\quad \E \left[ \int_0^T  L(x_t, u_t)dt\right]
\end{equation}
subject to the constraint $(x_t, u_t)\in X\times U$ with the initial condition $x_0\in X$ and the terminal condition $x|_{t=T}=x_T\in X$ with probability $1$. Here, in general, $x_0$ is an $\mathcal{F}_0$-measurable random variable and $x_T$ is an $\mathcal{F}_T$-measurable random variable. An optimal path $(x_t^*, u^*_t)$ is an $\mathcal{F}_t$-adapted process such that \eqref{DS_1} or \eqref{Ctn_1} is achieved.

In the case where the terminal value is not given and a terminal payoff function $S$ is considered, the problems can be reformulated as
\begin{equation}\label{DS_3}
\text{maximise}\quad \E\left[ \sum_{t=0}^T  L(x_t,u_t)+S(x_T)\right]
\end{equation}
and
\begin{equation}\label{Ctn_3}
\text{maximise}\quad \E \left[ \int_0^T  L(x_t, u_t)dt+S(x_T)\right]
\end{equation}

Further, when the future is discounted, the problem is described, in the discrete-time setting as
\begin{equation}\label{DS_2}
\text{maximise}\quad \E\left[ \sum_{t=0}^T \delta^t L(x_t,u_t)\right]
\end{equation}
and in the continuous-time setting as
\begin{equation}\label{Ctn_2}
\text{maximise}\quad \E \left[ \int_0^T  e^{-\rho t}L(x_t, u_t)dt\right].
\end{equation}

The stochastic analog of the turnpike theorem states that the optimal paths converge in a stochastic sense (say, convergence in distribution, convergence in probability or almost surely convergence) to a unique stochastic limit under standard assumptions of differentiability and concavity of the function $L$. The proof was mostly based on dynamic programming principle.

The most studied stochastic version of the middle turnpike theorem states that for any $\epsilon>0$, there exists a number $\theta>0$ such that the expected number of time periods for which $\{x_t^*\}$ spends outside the $\epsilon$-neighbourhood of $x^*$, i.e. $||x_t^*-x^*||>\epsilon$, is bounded from above by $\theta$, where $\theta$ is independent of $T$.



The uniformity assumption and the sensitivity assumption were mostly used to prove the stochastic version of turnpike theorems, as discussed in Joshi \cite{J2003}. The idea is first to construct a value-loss process, say $\{V_t:t\geq 0\}$ which has
the convenient property of being either a martingale (Majumdar and Zilcha \cite{MZ1987}), a submartingale (Brock and Majumdar \cite{BM1978} or a supermartingale (Joshi \cite{J1997}) thereby permitting a passage to the rich theory of martingales. Then the process $\{V_t:t\geq 0\}$ is assumed to be uniform bounded in expectation (i.e. the uniformity assumption) and have certain sensitivity property (i.e. the sensitivity assumption). The sensitivity assumption of $\{V_t:t\geq 0\}$ says that, for any $\epsilon>0$, if the optimal paths diverge from the stochastic limit by more than $\epsilon$, then the process $\{V_t\}$ records a value-loss of at least  $\eta(\epsilon)>0$. In other words, value-loss should be sensitive to a critical divergence of the optimal paths.

For general Markov chains the question on the existence of optimal value for discounted prices wwas addressed in \cite{Hol86}.
Turnpike theorems for Markov chains were proved by Shapiro \cite{Sha68} in discounted case and by Kolokoltsov \cite{Kol89} without discounting. In the latter paper the turnpike theorems on the strategies and on the distributions of states were distinguished.


For the financial portfolio problems, the turnpike theorems was proved, mainly under convexity conditions on the function $L$ directly,  e.g. by   Mossin \cite{M1968}, Leland \cite{L1972}, Hakansson \cite{H1974}, Hakansson \cite{H1974}, Sumit Joshi \cite{J1997}, Cox and Huang \cite{CH1992} and Huang and Zariphopoulou \cite{HZ1999}, Dybvig, Rogers and Back Kerry \cite{DRK1999}. In the discrete-time setting, Mossin \cite{M1968} showed the convergence of optimal strategies under the assumption that $-\frac{L'(x)}{L''(x)}=ax+b$, where $L'$ and $L''$ is the first and second order derivative of $L$, respectively. Leland \cite{L1972} generalised Mossin's result by assuming that $-\frac{L'(x)}{L''(x)}=ax+h(b)$
where the function $h$ is uniformly bounded. Hakansson \cite{H1974} proved convergence of optimal strategies by making assumptions slightly more general $\frac{(x-a)^{1-b}}{1-b}\leq L(x)\leq \frac{(x-a)^{1-b}}{1-b}, \forall x\geq a $. Huberman and Ross \cite{HR1983} showed that convergence occurs for the functions $L$ that are bounded from below and for some $c\in(0,1)$, $\lim_{x\to \infty}\frac{- L''(x)x}{L'(x)}=c$.

In the continuous-time setting, Cox and Huang \cite {CH1992} demonstrate that the portfolio turnpike property holds if the utility function $L$ is strictly concave with $\lim_{x\to \infty}L'(x)=0$ and there exists $A_1,A_2, b>0$ and $K>0$ such that
$|(L')^{-1}(x)-A_1x^{-\frac{1}{b}}|\leq A_2 x^{-a}$, for all $x\leq K$ and some $a\in [0, \frac{1}{b})$, where $(L')^{-1}$ denotes the inverse of $L'$. Huang and Zariphopoulou \cite{HZ1999} gave another sufficient condition for the turnpike property that the strictly concave function $L$ satisfies the condition $\lim_{x\to \infty}\frac{L'(x)}{x^{r-1}}=k$ for $0<r<1$ and $k>0$.

In optimal growth theory, mainly under the uniformity assumption and the sensitivity assumption, in the undiscounted case, Brock and Mirman \cite{BM1973}, Dana \cite{D1974}, Evstigneev \cite{E1974, E1976, AE1987} and Mirman and Zilcha \cite{MZ1977} showed that all  optimal paths converge in an appropriate topology to the optimal stationary path. In the discounted context, Brock and Mirman \cite{BM1972} and Mirman and Zilcha \cite{MZ1975} proved the convergence in distribution of optimal paths from distinct initial stocks to the stationary optimal paths. Fleming, Sethi and Soner \cite{FSS1987} studied production planning problem where the demand is described as a Markov chain and showed a turnpike theorem for the Markov control.

Some research considered a so-called {\it non-stationary} setting. The stochastic environment can be non-stationary, namely, $\Omega:=X_{t=0}^\infty \Omega_t$  and is considered as the set of all sequences $\omega:=\{\omega_t: \omega_t\in \Omega_t, t\geq 0\}$, where $\omega_t$ is a possible state of the environment at time $t$. If the function $L$ depends additionally on time $t$, it is called time-varying or non-stationary preference. In non-stationary models, the convergence in probability (Brock and Majumdar \cite{BM1978}) and the stronger property of almost sure convergence (Chang \cite{C1982}, F\"ollmer and Majumda \cite {FM1978}, Joshi \cite{J1997}, Majumdar and Zilcha \cite{MZ1987}, Joshi  \cite{J2003}) were proved.  Denardo and Rothblum \cite{DR2006} considered an exponential utility function $L$ for a Markov decision process and proved that if the corresponding transition rate matrix is transient, there exits an optimal stationary control path $u^*$.


In game theory, turnpike theorems for long-time two-player zero-sum  stochastic games on a finite state space were proved by e.g. Kolokoltsov \cite{K1992}, Kolokoltsov and Malafeyev\cite{KM2010}, by utilising the Bellman (or Shapley) operator  and imposing conditions on the transition probability. Related research on  long-time two-player zero-sum  stochastic games can be founded e.g. in Sorin \cite{S2005}, Hern\'andez-Lerma and  Lasserre \cite{HL2001} , Guo and Hern‡ndez-Lerma \cite{GH2005}.

Haurie and Delft \cite{HD1991}  considered a dynamics games  with piecewise deterministic path and random jumps. In the discounted setting,  under the assumption of convexity and piecewise uniqueness, the global asymptotic was proved if the transition probability $q_{ij}$ are constant and the local asymptotic was proved  if $q_{ij}$ are state and control dependent.

For zero-sum Markov games, Ja\'skiewicz  \cite{J2002},  Vega-Amaya \cite{VA2003}
worked with  strongly continuous transition probabilities, i.e. the transition probability $q(D|x, u, v )$ is continuous in actions $u$ and $v$, for any fixed states $x$ and Borel subset $D$ of the state space $X$. It was proved that the optimality equation has a solution, there exists a value of the game and both players possess optimal stationary strategies.  In the proof, Ja\'skiewicz  \cite{J2002} analyzes auxiliary perturbed models, whereas Vega-Amaya \cite{VA2003} makes use of a fixed point theorem.

In the work of Ja\'skiewicz and Nowak \cite{JN2006} and K\"uenle \cite{K2007}, zero-sum Markov games within the similar setting where the transition times occur at integer time moments were studied. Ja\'skiewicz and Nowak \cite{JN2006} prove the similar results  by applying Fatou's Lemma for weakly convergent measures. K\"uenle \cite{K2007}  introduced certain contraction operators which lead to a parameterized family of functional equations. Making use of some continuity and monotonicity properties of the solutions to these equations (with respect to the parameter) he obtains a lower semicontinuous solution of the optimality equation. Then he showed that the maximzing player has an $\epsilon$-optimal stationary strategy, and the minimizing player has an optimal stationary strategy.

Ja\'skiewicz \cite{J2009} studied a zero-sum ergodic semi-Markov games, in a general Borel state space, with weakly continuous (Feller) transition probabilities and lower semicontinuous payoff function $L$. Under the Banach fixed point theorem and an ergodicity assumption of an embedded Markov chain, it was proved that in this game, one player possesses an $\epsilon$-optimal stationary strategy and the other has an optimal stationary strategy.

Beyond the uniqueness assumption on the turnpike, it is worth noting that  if there are several optimal stationary paths, there exists a state transition of the optimal paths from one optimal stationary path to another, see e.g. Mamedov \cite{MM2003}, Alain Rapaprot and Pierre Cartigny \cite{RC2004}, Kamihigashi and Roy \cite{KR2006}.

\section{Basic setting for a zero-sum Markov game}\label{modelling}

In this section, we shall set up a model for a discrete-time two-person zero-sum Markov game on a general Borel state space. A Bellman (or Shapley)operator is defined for this game and some basic properties of this operator are discussed.

A discrete-time two-person zero-sum $T$-step, $T\in \N$, Markov game is a tuple $(X,\pi, U, V, P, g, S)$ of the following meanings:

$\bullet$ $X$ is a Borel subset of a general complete metric state space endowed with its Borel $\sigma$-field $\mathcal{B}(X)$ and a Borel $\sigma$-finite reference measure $\pi$ on it. One can always keep in mind the sets $\{1, \cdots, n\}$ and $\N$ with arbitrary discrete measure or a subset of $\R^d$ with Lebesgue measure as basic examples of $X$ with its reference measure $\pi$.

$\bullet$ $U$ and $V$ are compact spaces of the admissible strategies (i.e., control parameters) for player 1 and player 2, respectively.

$\bullet$ $P=\{ P(u,v,x,\cdot), x\in X, u\in U, v\in V \}$ is the family of transition probability kernels, i.e. $P(u,v,x,A)$, $A\in\mathcal{B}(X)$, is the transition probability from state $x$ to  $A$ if the two players choose strategies $u$ and $v$ respectively.  We assume that the transition probability $P(u,v,x, \cdot)$ has a density $\rho( u,v,x,y)$ with respect to the reference measure $\pi (\cdot)$.

$\bullet$ $g:U\times V\times X \times X\mapsto \R$ is a transition cost function for player 1, namely, $g(u,v,x,y)$ is the cost of the player 1 (i.e., the income of the player 2) from the transition from $x$ to $y$ under controls $u$ and $v$. If the cost $g(u,v,x,y)$ is negative, this means that the player 2 pays $|g(u,v,x,y)|$ to player 1. The function $g$ is assumed to be continuous and also bounded by certain constant $G>0$, uniformly in $u, v, x, y$, i.e.
$$|g(u,v,x,y)|\leq G, \quad \forall u,v,x,y.$$

$\bullet$ $S:X\mapsto \R$ is a final cost function depending on the final position. The function $S$ is assumed to be bounded and continuous. We denote the set of bounded and continuous functions on $X$ by $C_b(X)$.

The dynamic $T$-step game with a given initial position $x\in X$ is played as follows. At the first step, starting from the given initial state $x$, the players choose sequentially certain strategies $v\in V$ and then $u\in U$ and  the game moves to the next state $y\in X$ with probability $\rho(u,v,x,y)\pi(dy)$. Then the player 1 pays $g(u,v,x,y)$ to player 2 if $g(u,v,x,y)$ is positive; otherwise, the player 2 pays $|g(u,v,x,y)|$ to player 1. The second step is played analogously starting from the new position $y$. After $T$ steps, the player 2 receives additional payment $S(k)$ depending on the finial position $k\in X$. In this multi-step game, the player 1 aims to minimise the payoff to player 2.

We define the Bellman operator $B: C_b(X)\mapsto C_b(X)$ by
\begin{equation}
(BS)(x):=\inf_u \sup_v \left[\int_X \big(g(u,v,x,y)+S(y)\big)\rho(u,v,x,y)\pi(dy)\right].
\end{equation}
Then, according to the dynamic programming principle, the solution to the $T$-step Markov game is given by the iterations $B^TS$. We would like to mention that, since the transition cost function $g$ and the final cost function $S$ are assumed to be bounded and continuous, the $\inf\sup$ in the definition (\ref{Bellman operator2}) can be replaced by $\min\max$. Namely, in our setting, the Bellman operator $B$ can be written as
\begin{equation}\label{Bellman operator2}
(BS)(x)=\min_u \max_v \left[\int_X \big(g(u,v,x,y)+S(y)\big)\rho(u,v,x,y)\pi(dy)\right].
\end{equation}

\begin{Remark}
It is often assumed that the game defined by (\ref{Bellman operator2}) has a value, which means that $\inf\sup$ are interchangeable for all $S\in C_b(X)$ such that one also has
\[
	(BS)(x)= \sup_v \inf_u \left[\int_X \big(g(u,v,x,y)+S(y)\big)\rho(u,v,x,y)\pi(dy)\right].
\]
But we shall not make or use this assumption. Hence, in our framework, $B^T$ describes the payoff for the player 1 calculated for the worst scenario dynamics (sometimes called robust control).
\end{Remark}

It is clear that, for any constant $a\in\R$ and any function $S\in C_b(X)$, the Bellman operator satisfies the following {\it homogeneity property}
\begin{equation}\label{linearity2}
B(a+S)=a+BS.
\end{equation}
Moreover, the Bellman operator $B$ is {\it non-expansive} in the usual sup-norm, i.e.
\begin{equation} \label{non-expensive2}
||BS_1-BS_2||\leq ||S_1-S_2||, \quad \forall S_1,S_2 \in C_b(X)
\end{equation}
since the inequality
\begin{equation*}
\begin{split}
	||BS_1-BS_2||=&\sup_{x\in X} |(BS_1-BS_2)(x)|\\
	\leq & \sup_{x\in X} \sup_{u,v} \int _X\left | (S_1(y)-S_2(y))\right |\rho(u,v,x,y)\pi(dy)\\
    \leq & ||S_1-S_2||
\end{split}
\end{equation*}
holds for any $S_1, S_2\in C_b(X)$.

With the homogeneity property (\ref{linearity2}) and the non-expansive property (\ref{non-expensive2}) of the Bellman operator $B$, we can define the quotient space $\tilde{C_b}(X)$ of $C_b(X)$ with respect to constant functions. Namely, the equivalence classes are defined by
$$[S]:=\{S+a, \forall a\in \R: S\in C_b(X)  \}.$$
Let $\Pi:C_b(X)\mapsto \tilde{C_b}(X)$ be the natural projection. The quotient norm on $\tilde{C_b}(X)$ is defined by the formula
\begin{equation}\label{norm}
||\Pi(f)||=\inf_{a\in \R}||f+a||=\frac{1}{2}\left(\sup_{x\in X}f(x)-\inf_{x\in X}f(x)\right).
\end{equation}
It is clear that $\Pi$ has a unique isometric right inverse mapping $\Phi: \tilde{C_b}(X) \mapsto C_b(X)$ so that $\Pi \circ \Phi$ is an identity in $\tilde{C_b}(X)$. The image $\Phi(\tilde{C_b}(X))$ consists of all functions $S\in C_b(X)$ such that $$\sup_{x\in X}S(x)=-\inf_{x\in X}S(x).$$ Thus one can identify $\tilde{S}\in \tilde{C_b}(X)$ with its image $S=\Phi(\tilde{S})\in C_b(X)$. By the properties (\ref{linearity2}) and (\ref{non-expensive2}), the continuous quotient map $\tilde{B}: \tilde{C_b}(X) \mapsto \tilde{C_b}(X)$ is well defined.

\section{Main results} \label{S3}

In this section, we shall present the turnpike theorems for the game specified above. First, the uniqueness of the average cost and the stationary strategy are proved. Then, the turnpike theorems on the set of strategies and on the state space are proved.

To state the main results, we need some additional properties of the transition probabilities:
\begin{equation}\label{prop2}
\begin{split}
\hspace{-.5cm}\exists \delta>0: \exists A\in \mathcal{B}(X)\, \text{with}\, &\pi(A)>0: \forall x\in X,  \forall u\in U, \forall v\in V, \forall y\in A,\\[.5em]
&\rho(u,v,x,y) \geq \delta.
\end{split}
\end{equation}

First, we show that  the operator $\tilde{B}$ is a contraction in the following
\begin{Lemma}\label{contraction_B}
If the property (\ref{prop2}) holds, then
\begin{equation}\label{contractionproperty}
||\tilde{B}(\tilde{S_1})-\tilde{B}(\tilde{S_2})||\leq (1-\delta\pi(A))||\tilde{S_1}-\tilde{S_2}|| \quad \forall \tilde{S_1}, \tilde{S_2} \in \tilde{C_b}(X).
\end{equation}
\end{Lemma}

\proof
For arbitrary $S_1, S_2\in C_b(X)$ and arbitrary states $x,z\in X$,
\begin{equation*}
\begin{split}
&(BS_1-BS_2)(x)-(BS_1-BS_2)(z)\\
\leq &\max_{u,v} \int_X(S_1(y)-S_2(y))\rho(u,v,x,y)\pi(dy)\\
&\hspace{5em}-\min_{u,v}\int_X(S_1(y)-S_2(y))\rho(u,v,z,y)\pi(dy).\\
\end{split}
\end{equation*}
Denote
\begin{equation*}
\begin{split}
	(u_1, v_1)=\arg \max_{u,v} \int_X(S_1(y)-S_2(y))\rho(u,v,x,y)\pi(dy)
\end{split}
\end{equation*}
and
\begin{equation*}
\begin{split}
	(u_2, v_2)=\arg \min_{u,v} \int_X(S_1(y)-S_2(y))\rho(u,v,z,y)\pi(dy).
\end{split}
\end{equation*}
Then, we have
\begin{equation}\label{basic}
\begin{split}
&(BS_1-BS_2)(x)-(BS_1-BS_2)(z)\\
\leq& \int_X(S_1(y)-S_2(y))\Big(\rho(u_1,v_1,x,y)-\rho(u_2,v_2,z,y)\Big)\pi(dy)\\
\leq & ||S_1-S_2|| \int_X|\rho(u_1,v_1,x,y)-\rho(u_2,v_2,z,y)|\pi(dy).
\end{split}
\end{equation}
By property (\ref{prop2}), for all  $y\in A$,
$$\rho (u_1,v_1,x,y)\geq \delta, \quad \text{and}\quad \rho(u_2,v_2,z,y)\geq \delta.$$
Then
\begin{equation*}
\begin{split}
&\int_X|\rho(u_1,v_1,x,y)-\rho(u_2,v_2,z,y)|\pi(dy)\\
\leq& \int_{X\setminus A} \Big(\rho(u_1,v_1,x,y)+ \rho(u_2,v_2,z,y)\Big)\pi(dy)\\
&\hspace{5em}+\int_{A} \Big|\rho(u_1,v_1,x,y)- \rho(u_2,v_2,z,y)\Big|\pi(dy)\\
= &1-\int_{A} \rho(u_1,v_1,x,y)\pi(dy)+1-\int_{A}  \rho(u_2,v_2,z,y)\pi(dy)\\
&\hspace{5em}+\int_{A} \big|\rho(u_1,v_1,x,y)- \rho(u_2,v_2,z,y)\big|\pi(dy).\\
\end{split}
\end{equation*}
Denote by $A^+$ the set of $y\in A$ such that
$$ \rho(u_1,v_1,x,y)> \rho(u_2,v_2, z,y)$$
and
by $A^-$ the set of $y\in A$ such that
$$ \rho(u_1,v_1,x,y)\leq \rho(u_2,v_2,z,y).$$
Clearly, $A^+\cap A^-=\varnothing$ and $A^+\cup A^-=A$, so
\begin{equation*}
\begin{split}
&\int_X|\rho(u_1,v_1,x,y)-\rho(u_2,v_2,z,y)|\pi(dy)\\
\leq& 2-\int_{A^+\cup A^-}\rho(u_1,v_1, x,y)\pi(dy)-\int_{A^+\cup A^-}\rho(u_2,v_2, z,y)\pi(dy)\\
&\hspace{5em}+\int_{A^+}(\rho(u_1,v_1, x,y) - \rho(u_2,v_2, z,y))\pi(dy)\\
&\hspace{5em}+\int_{A^-}(\rho(u_2,v_2, z,y) - \rho(u_1,v_1,x,y))\pi(dy)\\
=&2-2\int_{A^+}\rho(u_2,v_2, z,y)\pi(dy)-2\int_{A^-}\rho(u_1,v_1, x,y)\pi(dy)\\
\leq &2-2\delta \pi(A^+)-2\delta \pi(A^-)\\
=&2-2\delta \pi(A).
\end{split}
\end{equation*}
Then, together with (\ref{basic}), we get
$$(BS_1-BS_2)(x)-(BS_1-BS_2)(z)\leq2(1-\delta\pi(A))||S_1-S_2||.$$
Now by the definition of the quotient norm (\ref{norm})
$$||\Pi(BS_1-BS_2)||\leq (1-\delta\pi(A))||S_1-S_2||.$$
So we obtain
$$||\tilde{B}\tilde{S_1}-\tilde{B}\tilde{S_2}||\leq (1-\delta\pi(A))||\tilde{S_1}-\tilde{S_2}||.$$
\qed

\begin{Theorem} [On the average cost]
If (\ref{prop2}) holds, then there exists a unique $\lambda \in \R$ and a function $S^*\in C_b(X)$, unique up to equivalence (i.e. $\tilde{S}^*=\Pi(S^*)$ is unique), such that
\begin{equation}\label{turnpike2}
B(S^*)=\lambda +S^*.
\end{equation}
Moreover, for every $S\in C_b(X)$,
\begin{equation} \label{estimates2}
||B^TS-T\lambda||\leq ||S^*||+||S^*-S||
\end{equation}
\begin{equation}\label{AC2}
\lim_{T\to \infty} \frac{B^TS}{T}=\lambda
\end{equation}
and, for every $\tilde{S} \in \tilde{C_b}(X)$,
\begin{equation}\label{convergence}
||\tilde{B}^T\tilde{S}-\tilde{S}^*||\leq (1-\delta\pi(A))^T ||\tilde{S}-\tilde{S}^*||.
\end{equation}
 \end{Theorem}

\proof
First, by  the contraction property (\ref{contractionproperty}) of the operator $\tilde{B}$ on $\tilde{C_b}(X)$ and the contraction principle (i.e. Banach fixed point theorem), the operator $\tilde{B}$ has a unique fixed point in $\tilde{C_b}(X)$. Namely, there exists a unique $\tilde{S}^*\in \tilde{C_b}(X)$ such that
\begin{equation}\label{fixpoint}
\tilde{B}\tilde{S}^*=\tilde{S}^*.
\end{equation}
Consequently, there exists a $\lambda\in \R$ and a function $S^*\in C_b(X)$, unique up to equivalence such that  (\ref{turnpike2}) holds.

Next, we shall prove (\ref{estimates2}). On one hand, by (\ref{turnpike2}) and (\ref{linearity2})
\begin{equation*}
\begin{split}
B^TS^*=B^{T-1}(S^*+\lambda)
=B^{T-2}(S^*+2\lambda)
=\cdots
=T\lambda+S^*
\end{split}
\end{equation*}
and together with the triangle inequality, for any $S\in C_b(X)$, we get
\begin{equation} \label{ineq21}
||B^TS-B^TS^*||=||B^TS-T\lambda-S^*||\geq ||B^TS-T\lambda||-||S^*||.
\end{equation}
On the other hand,  by (\ref{non-expensive2}), we have,
\begin{equation} \label{ineq11}
||B^TS-B^TS^*||\leq||B^{T-1}S-B^{T-1}S^*||\leq \cdots \leq||S-S^*||.
\end{equation}
By combining  (\ref{ineq21}) and (\ref{ineq11}), we have
$$ ||B^TS-T\lambda||-||S^*||\leq ||S-S^*||$$
which clearly gives the inequality (\ref{estimates2}).

Equation (\ref{AC2}) is implied by (\ref{estimates2}) by passing to the limit $T\to \infty$,
$$0 \leq \lim_{T\to \infty} \frac{||B^TS-T\lambda||}{T}\leq \lim_{T\to \infty} \frac{||S^*||+ ||S-S^*||}{T}=0$$
as the functions $S$ and $S^*$ are all bounded functions. Equation (\ref{AC2}), in turn, implies the uniqueness of $\lambda$ in (\ref{turnpike2}).

Further, for every $ \tilde{S} \in \tilde{C_b}(X)$, by (\ref{fixpoint}) we have
\begin{equation*}
\begin{split}
||\tilde{B}^T\tilde{S}-\tilde{S}^*|=&||\tilde{B}(\tilde{B}^{T-1}\tilde{S})-\tilde{B}
(\tilde{B}^{T-1}\tilde{S}^*)||\\
\leq&(1- \delta\pi(A))||\tilde{B}^{T-1}\tilde{S}-\tilde{B}^{T-1}\tilde{S}^*||\\
\leq&\cdots\\
\leq& (1- \delta\pi(A))^T||\tilde{S}-\tilde{S}^*||
\end{split}
\end{equation*}
which completes the proof.
\qed

\begin{Remark}
The unique $\lambda$ in (\ref{turnpike2}) can be interpreted as the average cost of player 1 from the long time Markov game.
\end{Remark}

Now we are ready for the first turnpike theorem on the set of optimal strategies.  Here we introduce some notations. Denote by $E(x,S)=\{(u(x),v(x))\}$ the set of pairs of optimal (equilibrium) strategies $(u(x),v(x))$ of (\ref{Bellman operator2}) for any $S\in C_b(X)$. Let $E(S)=\cup_{x\in X}\{E(x,S)\}$ be the set of pairs of optimal (equilibrium) strategies with final cost function $S\in C_b(X)$.  In particular, the family $E(S^*)$ contains the sets of the stationary strategies in the game and $E(B^{T-t}S)$ is that the family of the sets of optimal (equilibrium) strategies at the step $t$ of a $T$-step game with final cost function $S\in C_b(X)$.

\begin{Theorem}[Early turnpikes on the set of strategies] \label{TH1}
Let (\ref{prop2}) hold. Assume the set-valued function $E(\cdot)$ is continuous at least at the point $S=S^*$. Then for an arbitrary neighbourhood $U(E(S^*))$ of the set $E(S^*)$, there exists an $M\in N$ such that if $T>M$, then for any $S\in C_b(X)$
$$E(B^{T-t}S)\subset U(E(S^*))$$
for all $t<T-M$.
\end{Theorem}

\begin{Remark}
The continuity of the set-valued function $E(\cdot)$ is in the sense of Hausdorff metric. The property that $E(\cdot)$ is continuous at the point $S=S^*$ has been used and discussed in many places. For example, in \cite{HMC2006}, this property is called the sensitivity property for a control law.
\end{Remark}

\proof
By (\ref{convergence}), we have for any fixed $T$
$$\tilde{B}^{T-t} \tilde{S} \rightarrow \tilde{S}^* \quad \text{for large enough} \,\, T-t.$$
By the assumption that $E(S)$ depends continuously on $S$ around $S^*$, we have
$$E(\tilde{B}^{T-t} \tilde{S}) \rightarrow E(\tilde{S}^* )\quad \text{for large enough} \,\, T-t.$$
So there exists an $M\in \N$, such that
$$E(\tilde{B}^{T-t} \tilde{S})\subset U(E(\tilde{S}^*)),\quad t<T-M$$
which completes the proof.
\qed

\begin{Remark}
This theorem states that , in a sufficiently large $T$-step game with an arbitrary final cost function $S\in C_b(X)$,  the optimal strategies at steps $t<T-M$ must be close to the stationary strategies $E(S^*)$ and diverge away from the turnpike only near to the end of the game. Thus this is a stochastic competitive  control version of { \it the early turnpike} as defined in \cite{Mc1976}.
\end{Remark}

Further, with additional assumptions on the set of optimal strategies and on the equilibrium transitions, we shall show that, in a long-time game, if the game is carried out with optimal strategies, then, apart from short initial and final periods of time, the distributions of its position on the state space must be very close to the stationary distribution. Let $\|\cdot\|_{TV}$  denote the total variation norm of measures.

\begin{Theorem} [Middle turnpike on the (distribution of) state space]
Let (\ref{prop2}) hold. Assume the set-valued function $E(S)$ is continuous at least at the point $S=S^*$ and, for each state $x\in X$, the set $E(x,S^*)$ contains only one pair of optimal strategies $(u^*(x), v^*(x))$. Let $Q^*=\{q^*(x)\}_{x\in X}$ denote the stationary distribution for the stationary Markov chain defined on the state space $X$ by the unique pair of optimal strategies $(u^*(x), v^*(x))$. Assume also that the equilibrium transition probabilities depend locally Lipschitz continuous on $S$ around $S^*$, i.e. that there exists a constant $k$ such that,
\begin{equation}\label{plip}
 ||P(u(x),v(x),x,\cdot)-P(u^*(x),v^*(x), x,\cdot)||_{TV}\leq k||S-S^*||,\quad \forall x\in X
\end{equation}
where $(u,v)\in E(x,S)$ and $S$ is sufficiently close to $S^*$. Then for each $\epsilon>0$, there exists an $M\in N$ such that for each $T$-step game, $T>2M$, with terminal cost $S\in C_b(X)$ of  player 1, we have
$$||Q(t)-Q^*||_{TV}<\epsilon$$
for all $t\in [M, T-M]$, where $Q(t)=\{q(x,t)\}_{x\in X}$ and $q(x,t)$ is the probability that the process is in a state $x\in X$ at time $t$ if the game is carried out with the optimal strategies.
\end{Theorem}
\begin{Remark}
In other words, $q^*(x)$ is the mean amount of time that each sufficiently long game with optimal strategies spends in position $x$.
\end{Remark}

\proof
From (\ref{plip}), together with Lemma \ref{contraction_B}, we have
\begin{equation*}
 ||P_t(x,\cdot)-P^*(x,\cdot)||_{TV}\leq k\epsilon^{T-t}||S-S^*||
\end{equation*}
where $\epsilon=1-\delta \pi(A)$, $P_t(x,\cdot)$ is the transition probability arising from optimal strategies $(u(x),v(x))$ at time $t\leq T$ with terminal payoff $S$ and $P^*(x,\cdot)=P(u^*(x),v^*(x), x,\cdot)$.

Let $\Psi_t$ and $\Psi^*$ denote linear operators on the set of measures
$$(\Psi_t\mu)(dx):=\int\mu(dz)P_t(z,dx)\,\quad \text{and}\,\quad (\Psi^*\mu)(dx):=\int\mu(dz)P^*(z,dx)$$
so that
\begin{equation*}
\begin{split}
||\Psi_t-\Psi^*||=&\sup_{||\mu||=1} \int\mu(dz)\left(P_t(z,dx)-P^*(z,dx) \right)\\
    &\leq k\epsilon^{T-t}||S-S^*||.
\end{split}
\end{equation*}
The state distribution at time $t$ can be estimated by
\begin{equation*}
\begin{split}
Q(t)=&\Psi_t\cdots\Psi_1 Q(0)\\
	=&(\Psi^*+\xi_t)\cdots (\Psi^*+\xi_1)Q(0)\\
\end{split}
\end{equation*}
with $\xi_j\leq k\epsilon^{T-j}$, where $Q(0)$ is the initial state distribution. Since
\begin{equation*}
\begin{split}
    &||\Psi_t\cdots\Psi_1-(\Psi^*)^t||\\
	=&|| (\Psi_t-\Psi^*) \Psi_{t-1}\cdots\Psi_1 +\Psi^*(\Psi_{t-1}-\Psi^*) \Psi_{t-2}\cdots\Psi_1 \\
    &\hspace{3em}+\cdots+ (\Psi^*)^{t-1}(\Psi_1-\Psi^*)||\\
\leq & k\left[\epsilon^{(T-t)}+\epsilon^{(T-t+1)}+\cdots+\epsilon^{(T-1)}\right]\\
\leq & k \epsilon^{(T-t)} (1+\epsilon+\epsilon^2+\cdots)\\
=&\frac{k}{1-\epsilon}\epsilon^{(T-t)} \leq \frac{k}{1-\epsilon}\epsilon^{M_1}
\end{split}
\end{equation*}
for $t<T-M_1$, $M_1\in \N$, we have
$$||Q(t)-(\Psi^*)^tQ(0)||_{TV}\leq \epsilon_1$$
for $t<T-M_1$, where $\epsilon_1$ can be made arbitrary small by choosing large enough $M_1\in \N$.

By the theorem on the convergence of probability distribution in homogeneous Markov chains to a stationary distribution, it follows that for each $\epsilon_2>0$  there exists  an $M_2\in\N$ such that
$$||(\Psi^*)^tQ(0)-Q^*||_{TV}\leq \epsilon_2, \quad  t>M_2.$$
Let  $M=\max\{M_1,M_2\}$. If $T>2M$, the intersection of the sets $[1,T-M_1]$ and $[M_2, \infty)$ is non-empty. Then we have
$$||Q(t)-Q^*||_{TV}\leq\epsilon_1+ \epsilon_2\quad \forall t\in[M, T-M]$$
which completes the proof. \qed

\section{Concluding remarks}\label{conclusion}

We have considered a two-person zero-sum Markov game with a general Borel state space and the modelling is formulated in the discrete time setting. With the Bellman operator defined by (\ref{Bellman operator2}), we have proved the convergence of the long time average cost (reward) of the game and turnpike theorems on the set of optimal strategies and on the state space are proved. We would like to mention that these results still hold under the corresponding continuous time setting by discretising, that is by looking at a continuous time game as the limit of games with discrete times.

Notice that in our modelling, at each step, we can observe the exact position of this game. But it often happens that one has to take into account not only the exact position, but also, the state distribution (say, its variance of volatility). This leads to the turnpike properties for measure-valued Markov games, see \cite{K2010, K2012, KLY2012}. For example, the final cost function may be any continuous function on a measure space and the stochastic transition kernels may depend on the state distribution. We shall develop results under this more general setting in future work.


{\small
\begin{thebibliography}{99}

\bibitem{A1965} Atsumi, H. (1965). Neoclassical Growth and the efficient program of capital accumulation. {\it Review of Economic Studies}, 32(2), 127-36.

\bibitem{AB2007} Alvarez O., M. Bardi (2007). Ergodic problems in differential games. \textit{Advances in Dynamic Game Theory}, S. Jorgensen, M. Quincampoix, and T.L. Vincent, eds., pp 131-152, Ann. Internat. Soc. Dynam. Games, vol. 9, Birkhaeuser, Boston.

\bibitem{AB2010} Alvarez O.,  M. Bardi (2010). Ergodic Ergodicity, stablization, and singular perturbations of Bellman-Isaacs equations, \textit{Memoirs American Math. Soc.}, 77 pages.

\bibitem{AE1987} Arkin V. I., I.V. Evstigneer (1987). { \it Stochastic Models of Control and Economic Dynamics (Economic Theory, Econometrics, and Mathematical Economics)}. Academic Press.

\bibitem{AS1977} Araujo, A. Scheinkman, J.A. (1977). Smoothness, comparative dynamics, and the turnpike property. Econometrica, 45(3), 601-620.

\bibitem{B1982} Bewley T.F. (1982). An integration of equilibrium theory and turnpike theory. \textit{Journal of mathematical economics}, 20, 233-267.

\bibitem{B2007} Bewley, T.F.(2007). General equilibrium, Overlapping Generations Models and Optimal growth theory, Harvard University Press, Cambridge.

\bibitem{B2009} Bardi M. (2009). On differential games with long-time-average cost. \textit{Annals of the International Society of Dynamic Games}, 10(1), 1-16.

\bibitem{BC2007} Blot J., Crettez, B (2007). On the smoothness of optimal paths II: some local turnpike results. \textit{Dicisions Econ Finan}, 30, 137-150.

\bibitem{BM1972} Brock, W.A., Mirman, L.J. (1972). Optimal economic growth and uncertainty: the discounted case. \textit{Journal of Economic Theory}, 4, 479-513.

\bibitem{BM1973} Brock, W.A., Mirman, L.J. (1973). Optimal economic growth and uncertainty: the no discounting case. \textit{International Economic Review}, 14(3), 560-573.

\bibitem{BM1978} Brock, W.A., Majumdar, M. (1978). Global asymptotic stability results for multisector models of optimal growth under uncertainty when future utilities are discounted. \textit{Journal of Economic Theory}, 18, 225-243.

\bibitem{c1966} Cass D. (1966). Optimum Growth in an Aggregative Model of Capital Accumulation: A Turnpike Theorem. \textit{Econometrica}, 34(4), 833-850.

\bibitem{C1990}  Carlson D.A.(1990). The existence of catching-up optimal solutions for a class of infinite horizon optimal control problems with time delay. \textit{SIAM Journal on Control and Optimization}, 28, 402-422.

\bibitem{CJH1987} Carlson D.A., Jabrane A., Haurie A. (1987). Existence of overtaking solutions to infinite dimensional control problems on unbounded time intervals. \textit{SIAM Journal on Control and Optimization}, 25, 1517-1541.

\bibitem{CH1995} Carlson D.A., Haurie A.B.  (1996). A turnpike theory for infinite horizon open-loop differential games with decoupled controls. \textit{SIAM Journal on Control and Optimization}, 34(4), 1405-1419.

\bibitem{C1982} Chang, F.R. (1982). A note on the stochastic value loss assumption. \textit{Journal of Economic Theory}, 26, 164-170.

\bibitem{C1985} Coles, J.L. (1985). Equilibrium turnpike theory with constant returns to scale and possibly heterogeneous discount factors. {\it International Economic Review}, 26(3), 671-679.

\bibitem{CH1992} Cox J.C., C.-F. Huang (1992). A continuous-time portfolio turnpike theorem. {\it Journal of Economic Dynamics and Control}, 16, 491-507.

\bibitem{D1974} Dana, R.A.(1974). Evaluation of development programs in a stationary stochastic economy with bounded primary resources. In: Proceedings of the Warsaw Symposium on Mathematical Methods in Economics. North-Holland, Amsterdam, 179-205.

\bibitem{DMc1985} Dasgupta, S., Mckenzie, K.W. (1985). A note on comparative dynamics of stationary states. {\it Econ. Lett}, 18, 333-338.

\bibitem{DR2006} Denardo E.V.,  Rothblum U.G. (2006). A turnpike theorem for a risk-sensitive Markov decision process with stopping. {\it SIAM J. Control Optimal}, 45(2), 414-431.

\bibitem{D2005} Dindo P. (2005). A tractable evolutionary model for the Minority Game with asymmetric payoffs. \textit{Physica}, A (355), 110-118.

\bibitem{DRK1999} Dybvig P.H., Rogers L.C. G., Back Kerry (1999). Portfolio Turnpikes, the Review of Financial Studies, 12(1), 165-195.

\bibitem{E1974} Evstigneev I.V. (1974).  Optimal stochastic programs and their stimulating prices. {\it mathematical Models in Economics}, J. Los and M.W. Los (eds.), 219-252, North-Holland.

\bibitem{E1976} Evstigneev I.V. (1976).  Turnpike theorems in probabilistic models of economic dynamics. Mathematical Notes 19(2), 165-171, (translated from Matematicheskie Zametki, Vol. 19(2), 279-290, 1976).


\bibitem{DSS1958} Dorfman R., P. Samuelson and R. Solow (1958). \textit{Linear programming and economic analysis}. New York: McGraw-Hill.

\bibitem{FK1990} Fershtman C., M. I. Kamien (1990). Turnpike properties in a finite horizon differential game: dynamic duopoly with sticky prices. \textit{International Economic Review}, 31(1), 49-60.

\bibitem{FM1978} F\"ollmer, H., Majumdar, M. (1978). On the asymptotic behaviour of stochastic economic processes. {\it Journal of Mathematical Economics}, 5, 275–287.

\bibitem{FO1985} Feinstein C.D., Oren S.S. (1985). A "funnel" turnpike theorem for optimal growth problems with discounting. {\it Journal of Economic Dynamics and Control}, 9, 25-39.

\bibitem{FSS1987}  Fleming W., S.P. Sethi, Soner H.M. (1987). An optimal stochastic production planning problem with randomly fluctuating demand.{\it SIAM J. control Optim}, 25, 1494-1502.

\bibitem{G1967} Gale, D. (1967). On optimal development in a multi-sector economy. {\it Rev. Econ. Stud.}, 34, 1-18.

\bibitem{GH2005}  Guo X. and O. Hern\'andez-Lerma (2005). Zero-sum continuous-time Markov games with unbounded transition and discounted payoff rates. {\it Bernoulli}, 11, 1009-1029.

\bibitem{H1974} Hankansson, N. (1974).Convergence to isoelasitic utility and policy in multiperiod portfolio choice. {\it Journal of Financial economics}, 1, 201-224.

\bibitem{H1976} Haurie A. (1976). Optimal control on an infinite time horizon: the turnpike approach. {\it Journal of mathematical economics}, 3(1), 81-102.

\bibitem{HD1991} Haurie A.,  Delft C.V.(1991). Turnpike properties for a class of piecewise deterministic systems arising in manufacturing flow control.{\it  Annals of operations Research}, 29, 351-374.

\bibitem{HL2001}  Hernandez-Lerma O. and J. B. Lasserre (2001).  Zero-Sum Stochastic Games in Borel Spaces: Average Payoff Criteria, {\it SIAM J. Control Optim.}, 39, 1520-1539.

\bibitem{Hol86} U. D. Holzbaur. Fixed point theorems for discounted finite Markov decision processes.
{\it J. Math. Anal. Appl.} 116:2 (1986), 594-597.

\bibitem{HMC2006} Huang M.Y., R.P. Malhame, P. E. Caines (2006). Large population stochastic dynamic games: closed-loop Mckean-Vlasov systems and the Nash certainty equivalence principle. \textit{Communications in Information and Systems}, 6(3), 221-251.

\bibitem{HZ1999} Huang C.Fu., Zariphopoulou T. (1999). Turnpike behaviour of long-term investments. {\it Finance Stochastic}, 3, 15-34.

\bibitem{HR1983} Huberman G. and Ross S. (1983). Portfolio turnpike theorems, risk aversion and regularly varying functions. {\it Econometrica}, 51, 1345-1361.

\bibitem{J2002} Ja\'skiewicz A. (2002). Zero-sum semi-Markov games. {\it SIAM J. Control Optim}, 41, 723-739.

 \bibitem{J2009} Ja\'skiewicz A. (2009). Zero-sum ergodic semi-Markov games with weakly continuous transition probabilities. {\it J Optim Theory App}, 141, 321-347.

\bibitem{JN2006} Ja\'skiewicz A., Nowak, A.S. (2006). Zero-sum egodic stochastic games with Feller transition probabilities. {\it SIAM J. control Optimal}, 45, 773-789.

\bibitem{J1997}  Joshi S. (1997). Turnpike theorems in Nonconvex nonstationary Environments. {\it International econmic review}, 38(1), 225-248.

\bibitem{J2003} Joshi S. (2003). The stochastic turnpike property without uniformity in convex aggregate growth models. {\it Journal of economic Dynamics and control}, 27, 1289-1315.

\bibitem{KR2006} Takashi K., Roy S. (2006). Dynamic optimization with a non smooth, nonconvex technology: the case of a linear objective function. \textit{Economic Theory}, 29, 325-340.

\bibitem{KP2011}  Khan M. A. and Piazza A. (2011).  An overview of turnpike theory: towards the discounted deterministic case. {\it Adv. Math. Econ}, 14, 39-67.

\bibitem{Kol89} Kolokoltsov V.N. Turnpikes and infinite extremals in Markov decision processes. (Russian) {\it Mat. Zametki} 46:4 (1989),  118–120.

\bibitem{K1992} Kolokoltsov V.N. (1992). On linear, additive, and homogeneous operators in Idempotent Analysis.  \textit{Advances in Soviet Mathematics}  13, Idempotent Analysis, Ed. V.P.Maslov et S.N. Samborski,  87-101.

\bibitem{KM2010} Kolokoltsov V.N. and O.A. Malafeyev (2010). Understanding Game Theory: Introduction to the analysis of many agent systems of competition and cooperation. World Scientific.

\bibitem{K2010} Kolokoltsov V.N. (2010). Nonlinear Markov Processes and Kinetic Equations. Cambridge University Press, Cambridge.

\bibitem{K2012} Kolokoltsov V.N. (2012). Nonlinear Markov games on a finite state space (mean-field and binary interactions). to appear in {\it International Journal of Statistics and Probability}.

\bibitem{KLY2012} Kolokoltsov V.N. , J. Li and W. Yang (2012). Mean Field Games and Nonlinear Markov Processes. arXiv:1112.3744v


\bibitem{K2007} K\"uenle, H.U.(2007). On Markov games with average reward criterion and weakly continuous transition probabilities. {\it SIAM J. Control Optim}, 45, 2156-2168.

\bibitem{L1972} Leland, H. (1972). On turnpike portfolios, In G. Szego and K. Shell, (eds.), Mathematical Methods in Investment and Finance. Amsterdam: North-holland.

\bibitem{LL2000} Guerrero-Luchtenberg, C\'esar Leonardo (2000). A uniform neighborhood turnpike theorem and applications. {\it J. Nath. Econm.}, 34(3), 329-357.

\bibitem{LVM2006}  Le Van C.,  Morhaim L.(2006). On optimal growth models when the discount factor is near $1$ or equal to $1$. {\it International Journal of economic theory}, 2(1), 55-76.

\bibitem{M2002} Mamedov  M. A. (2002). Asymptotical stability of optimal paths in nonconvex problems, Research Report 02/02, University of Ballarat, February .

\bibitem{MM2003} Mamedov, M. A. (2003). A turnpike theorem for continuous-time control systems when the optimal stationary point is not unique.  {\it Abstract and Applied Analysis},  11, 631-650.


\bibitem{MM1982} Majumdar M., Mitra T. (1982). Intertemporal allocation with a Non-convex technology: the aggregative framework. \textit{Journal of Economic theory}, 27, 101-136.

\bibitem{MM1999} Marena, M., Montrucchio, L.(1999). Neighborhood turnpike theorem for continuous time optimization models. {\it J. Optim. Theory appl.}, 101, 651-676.

\bibitem{MN1982}  Majumdar and Nermuth (1982). Dynamic Optimisation in Non-Convex Models with Irreversible investment: Monotonicity and Turnpike Theorems. {\it Journal of Economics}, 42, 339-362

\bibitem{MR1973} Makarov V.L. and A. M. Rubinov (1973). Mathematical Theory of Economic Dynamics and Equilibria, Nauka, Moscow, English trans.: Springer-Verlag, New York.

\bibitem{MS1993} Maitra A. and Sudderth W. (1993). Borel stochastic games with limsup payoffs. {\it Ann. probab.}, 21, 861-885.

\bibitem{Mc1963} McKenzie L. W. (1963). The turnpike theorem of Morishima.  \textit{Review of Economic Studies}, 30, 169-176.

\bibitem{Mc1976} McKenzie L. W. (1976). Turnpike theory.\textit{Econometrica}, 44(5), 841-865.

\bibitem{Mc1977} McKenzie L. W. (1977).  A new route to the turnpike. In Henn, R. Moeachlin, O. (eds.) Mathematical Economics and Games Theory, Springer, New York

\bibitem{M1979}    Mckenzie, L. W. (1979). Optimal Economic growth and turnpike theorems. Disussion paper 79-1, uiversity of Rochester.

\bibitem{Mc1986} McKenzie L. W. (1986). Optimal economic growth, and turnpike theorems and comparative dynamics, In: Handbook of Mathematical Economics, 3, 1281-1355.

\bibitem{M1998} McKenzie, L.W.(1998). Turnpikes. {\it The American Economic Review, Papers and Proceedings of the Hundred and Tenth Annual meeting of the American Economic Association}. 88(2), 1-14.

\bibitem{M1995} Montrucchio L. (1995). A turnpike theorem for continuous-time optimal-control models.  {\it Journal of economic Dynamics and Control}, 19, 599-619.

\bibitem{M1995-2}  Montrucchio L. (1995). A new turnpike theorem for discounted programs.  {\it Econ. theory}, 5, 371-382.


\bibitem {MZ1975} Mirman, L.J., Zilcha, I. (1975). On optimal growth under uncertainty. {\it Journal of Economic Theory}, 11, 329-339.

\bibitem{MZ1977} Mirman, L.J., Zilcha, I. (1977). Characterizing optimal policies in a one-sector model of economic growth under uncertainty.  {\it Journal of Economic Theory}, 14, 389-401.

\bibitem{MZ1987} Majumdar, M., Zilcha, I. (1987). Optimal growth in a stochastic environment: some sensitivity and turnpike results.  {\it Journal of Economic Theory}, 43, 116-133.

\bibitem{M1961}  Morishima, Michio (1961). Proof of a Turnpike theorem: the no joint production case, \textit{Review of Economic Studies}, 28, 89-97.

\bibitem{M1968} Mossin, J.(1968). Optimal multiperiod policies.  {\it J. Business}, 31, 215-229.

\bibitem{P2000}  Park H.  (2000). Global asymptotic stability of a competitive equilibrium with recursive preferences.  {\it Economic theory}, 15(3), 565-584.

\bibitem{R1961} Radner, Roy (1961). Paths of Economic growth that are optimal with regard only to Final States. \textit{Review of Economic Studies}, 28, 98-104.


\bibitem{RC2004} Alain Rapaprot and Pierre Cartigny (2004). Turnpike theorems by a value function approach. \textit{SAIM: control, Optimisation and Calculus of Variations}, 10, 123-141.

\bibitem{SS1956} Samuelson, P.A., R.M. Solow (1956). A complete capital model involving heterogeneous capital goods.  {\it Q. J. Econ}, 27, 537-562.

\bibitem{Sha68} J. F. Shapiro. Turnpike Planning Horizon for a Markov Decision Model.
        {\it Management Science} 14 (1968), 292-300.

\bibitem{S2005} Sorin S. (2005). New approaches and recent advances in two-person zero-sum repeated games.  {\it Advances in Dynamic Games}, Nowak, Andrzej S.; Szajowski, Krzysztof (Eds.), 67-94,

\bibitem{S1976} Scheinkman J.(1976). On optimal steady states of $n$-sector growth models when utility is discounted. \textit{Journal of Economic theory}, 12, 11-30.

\bibitem{VA2003} Vega-Amaya, O. (2003). Zero-sum semi-Markov games. fixed-point solutions of the Shapley equation. \textit{SIAM J. control Optm}, 42, 1876-1894.

\bibitem{Y1984}  Yano  M.(1984). Competitive equilibria on turnpikes in a McKenzie Economy, I: A neighborhood turnpike theorem, {\it International Economic review}, 25(3) 695-717.

\bibitem{YK1992}  Yakovenko S.Yu., L.A. Kontorer (1992). Nonlinear semigroups and infinite horizon optimization, \textit{Advances in Soviet mathematics}, 13,167-210.

\bibitem{Z1995-1} Zaslavski  A.J. (1995). Optimal programs on infinite horizon 1. \textit{SIAM Journal on Control and Optimization}, 33, 1643-1660.

\bibitem{Z1995-2} Zaslavski  A.J. (1995). Optimal programs on infinite horizon 2, \textit{SIAM Journal on Control and Optimization}, 33, 1661-1686.

\bibitem{Z1999} Zaslavski  A.J. (1999). The turnpike property for dynamic discrete time zero-sum games. \textit{Abstract Applied Analysis}, 4(1), 21-48.

\bibitem{Z1998_2} Zaslavski  A.J. (1998). Turnpike theorem for convex infinite dimensional. \textit{Journal of convex analysis}, 5(2), 237-248.

\bibitem{Z2004} Zaslavski  A.J. (2004). Turnpike theorem for a class of discrete time optimal control problem. \textit{Fixed point theory and applications}, 5, 175-182.

\bibitem{Z2006} Zaslavski  A.J. (2006). Turnpike Properties in the Calculus of Variations and Optimal Control, Springer, New York.

\bibitem{Z2007} Zaslavski  A.J. (2007). Turnpike results for a discrete-time optimal control system arising in economic dynamics. \textit{ Nonlinear Analysis}, 67, 2024-2049.

\bibitem{Z2009} Zaslavski  A.J. (2009). Two turnpike results for a discrete-time optimal control system. \textit{ Nonlinear Analysis}, 71, 902-909.

\end{thebibliography}}
\end{document}